\documentclass[11pt]{amsart}
\usepackage{amssymb, latexsym}
\theoremstyle{plain}
\newtheorem{theorem}{Theorem}
\newtheorem{corollary}{Corollary}

\newtheorem{proposition}{Proposition}

\newtheorem*{2'}{Theorem 2'}
\newtheorem*{3'}{Theorem 3'}

\theoremstyle{remark}

\newtheorem*{Remark 1}{Remark 1}
\newtheorem*{Remark 2}{Remark 2}
\newtheorem*{Remark 3}{Remark 3}
\newtheorem*{Remark 4}{Remark 4}

\numberwithin{equation}{section}

\begin{document}

\title [Diffusive Search with Spatially Dependent Resetting]
{Diffusive Search with Spatially Dependent Resetting}

\author{Ross G. Pinsky}


\address{Department of Mathematics\\
Technion---Israel Institute of Technology\\
Haifa, 32000\\ Israel}
\email{ pinsky@math.technion.ac.il}

\urladdr{http://www.math.technion.ac.il/~pinsky/}

\subjclass[2010]{60J60} \keywords{random target, diffusive search, resetting,  optimization  }
\date{}

\begin{abstract}
We consider a stochastic search model with resetting for an unknown stationary target $a\in\mathbb{R}$ with known distribution $\mu$. The searcher begins at the origin and performs Brownian motion
with diffusion constant $D$. The searcher is also armed with an exponential clock with spatially dependent rate $r=r(\cdot)$, so that if
 it has  failed to locate the target by the time the clock rings, then its position is reset to the origin and it   continues its search anew from there. Denote the position of the searcher at time $t$  by $X(t)$.
Let $E_0^{(r)}$ denote expectations for the process $X(\cdot)$.
 The search ends at time $T_a=\inf\{t\ge0:X(t)=a\}$.  The expected time of the search is then $\int_{\mathbb{R}}(E_0^{(r)}T_a)\thinspace\mu(da)$.
Ideally, one would like to minimize this over all resetting rates $r$.
We obtain
quantitative growth
rates for $E_0^{(r)}T_a$  as a function of $a$ in terms of the asymptotic behavior of the rate function $r$, and also a rather precise
dichotomy on the  asymptotic behavior of the resetting function $r$ to determine whether $E_0^{(r)}T_a$ is finite or infinite.
We  show generically that if $r(x)$ is of the order $|x|^{2l}$, with $l>-1$, then
$\log E_0^{(r)}T_a$ is of the order $|a|^{l+1}$; in particular, the smaller the asymptotic size of $r$, the smaller the asymptotic growth
rate of $E_0^{(r)}T_a$. The asymptotic growth rate of $E_0^{(r)}T_a$ continues to decrease when $r(x)\sim \frac{D\lambda}{x^2}$
with $\lambda>1$; now the growth rate of $E_0^{(r)}T_a$ is more or less of the order $|a|^{\frac{1+\sqrt{1+8\lambda}}2}$.
Note that this exponent
increases to $\infty$ when $\lambda$ increases to $\infty$ and decreases to 2 when $\lambda$ decreases to 1. However, if $\lambda=1$, then
 $E_0^{(r)}T_a=\infty$, for $a\neq0$.
Our results  suggest that for many distributions $\mu$ supported on all of $\mathbb{R}$, a near optimal (or optimal) choice
of resetting function $r$ in order to minimize $\int_{\mathbb{R}^d}(E_0^{(r)}T_a)\mu(da)$ will be one
which decays quadratically as $\frac{D\lambda}{x^2}$ for some $\lambda>1$.
 We also give explicit, albeit  rather complicated, variational
formulas for $\inf_{r\gneqq0}\int_{\mathbb{R}^d}(E_0^{(r)}T_a)\mu(da)$.
For distributions $\mu$ with compact support, one should set $r=\infty$ off of the support. We also discuss
this case.
\end{abstract}

\maketitle
\section{Introduction and Statement of Results}\label{intro}
A number of recent papers have considered  a stochastic search  model for a stationary target  $a\in R^d$,
which might be random and have a known distribution attached to it, whereby a searcher sets off from a fixed point, say the origin, and performs  Brownian motion with diffusion constant $D$.
The searcher is also armed with   a (possibly space dependent) exponential resetting time, so that
if it has  failed to locate the target by the time the clock rings, then its position is reset to the origin and it continues  its search anew from there.
One may be interested in several statistics, the most important one being  the expected time to locate the target. (In dimension one, the target is considered ``located'' when the process
hits the point $a$, while in dimensions two and higher, one chooses an $\epsilon_0>0$ and the target is said to be ``located''
when the process hits the $\epsilon_0$-ball centered at $a$.) Without the resetting, this expected time is infinite.
When the resetting rate is constant,
the expected time to locate the target is finite.
  See, for example, \cite{EM1, EM2,EM3, EMM}. For related models with resetting, see \cite{G,MV1,MV2} as well as the   references in all of the above  articles.
Also see \cite{K} and \cite{P18} for a related problem without resetting, which is motivated by the above resetting problem.

Here is a more formal  mathematical definition of the model.
Let $a\in\mathbb{R}^d$ denote an unknown stationary target with known probability distribution $\mu\in\mathcal{P}(\mathbb{R}^d)$, the space of probability measures on $\mathbb{R}^d$.
The process $X(t)$ on $\mathbb{R}^d$ is defined as follows. Let $r=r(x)\gneqq0$ be a continuous function on $\mathbb{R}^d$.
This function will serve as the resetting rate.
The process starts from $0\in\mathbb{R}^d$ and  performs Brownian motion with diffusion constant $D$, until a random exponential clock rings.
The conditional probability that the clock has not rung by time $t>0$, given that the path up to time $t$ is $\{X(s),0\le s\le t\}$, is equal to
$\exp(-\int_0^t r\big(X(s)\big)ds)$.
When the clock rings, the process is instantaneously reset to its initial position 0, and continues its search  afresh with an independent resetting clock, and the above scenario is repeated, etc.
We define the process so that it is left continuous.
From the above description, it follows that  $X(\cdot)$ is a Markov process whose generator
$\mathcal{L}$  satisfies
\begin{equation}\label{generator}
\mathcal{L}u(x)=\frac D2\Delta u(x)+r(x)\big(u(0)-u(x)\big).
\end{equation}
(For more details on such types of constructions, see \cite{P09}.)

If $d=1$, let  $T_a=\inf\{t\ge0:X(t)=a\}$, while if $d\ge 2$, fix  $\epsilon_0>0$ and define $T_a=\inf\{t\ge0: |X(t)-a|\le \epsilon_0\}$.
Denote probabilities and expectations for the process starting at $x$ by $P^{(r)}_x$ and $E^{(r)}_x$ respectively.
Since the unknown target $a\in\mathbb{R}^d$ has   distribution $\mu$,
the expected search time is then given by  $\int_{\mathbb{R}^d}(E_0^{(r)}T_a)\mu(da)$.
Ideally, one would like to minimize this expression over all resetting rates $r$.

For most choices of $r$, it is not possible to write down a completely explicit expression for  $E_0^{(r)}T_a$.
If $r>0$ is constant, then one can calculate $E_0^{(r)}T_a$
explicitly in terms of appropriate Bessel functions \cite{EM3}. When $d=1$, this simplifies \cite{EM2} and one has
\begin{equation}\label{constantr}
E_0^{(r)}T_a=\frac{e^{\sqrt{\frac{2r}D}\thinspace |a|}-1}r,\ a\in\mathbb{R}.
\end{equation}
The only other case we've seen worked out explicitly in the literature is the case that $d=1$ and $r(x)$ is equal to 0 for $|x|<a_0$ and equal to a constant $r>0$
for $x\ge a_0$, where $a_0>0$ \cite{EM2}.

\it In this paper, we consider the one-dimensional case.\rm\ In Theorem \ref{rep}, for each $r\gneqq0$, we obtain an explicit formula  for $E_0^{(r)}T_a$ in terms of a positive function $\phi$ that solves $\frac D2\phi''-r\phi=0$.
It is intuitively clear that
$E_0^{(r)}T_a$ is not monotone in $r$. And indeed, one can see this explicitly  when $r$ is constant---it follows from \eqref{constantr} that
$E_0^{(r)}T_a$ approaches $\infty$ both as $r\to0$ and as $r\to\infty$. As will be seen, our explicit formula
for $E_0^{(r)}T_a$ in terms of the function $\phi$ is in fact the quotient of two functions, each of which is known
to be monotone in $r$.
This fact, along with   the explicit formula  for each of these two functions in the quotient, will allow us to obtain in Theorems \ref{comparison1} and \ref{comparison}  quantitative growth
rates for $E_0^{(r)}T_a$  as a function of $a$ in terms of the asymptotic behavior of the rate function $r$, and also in Theorem \ref{comparison} a rather precise
dichotomy on the asymptotic behavior of the resetting function $r$ which determines whether $E_0^{(r)}T_a$ is finite or infinite.
We also consider the case that the target distribution $\mu$ is compactly supported, in which case
it is advantageous for the searcher
to  be instantaneously reset
as soon as its position  has left the support of $\mu$.


We begin with a proposition which supplies us with several options for the above-mentioned function $\phi$.
\begin{proposition}\label{criticality}
Let $r\gneqq0$ be a continuous function on $\mathbb{R}$. Then there exist strictly positive functions $\{\phi_i\}_{i=1}^3$, all satisfying
\begin{equation}
\frac D2\phi''(x)-r(x)\phi(x)=0,\ x\in\mathbb{R},
\end{equation}
and such that
\begin{equation}\label{intcond}
\begin{aligned}
&\int_{-\infty}\phi_1^{-2}(x)dx=\infty,\ \  \int^{+\infty}\phi_1^{-2}(x)dx<\infty,\\
&\int_{-\infty}\phi_2^{-2}(x)dx<\infty,\ \  \int^{+\infty}\phi_2^{-2}(x)dx=\infty,\\
&\int_{-\infty}^{+\infty}\phi_3^{-2}(x)dx<\infty.
\end{aligned}
\end{equation}
Furthermore, if $r(x)$ is an even function, then $\phi_3(x)$ can be chosen to be even.
\end{proposition}
For $a>0$, let $u_{+,a}$ denote the solution to the equation
\begin{equation}\label{u+a}
\begin{aligned}
&\frac D2 u''-r(x)u=0, \ -\infty<x<a;\\
&u(a)=1;\\
&0\le u\le 1 \ \text{and}\ u\ \text{maximal}.
\end{aligned}
\end{equation}
(The condition  that $u$ be maximal means that any other  solution bounded above and below by 1 and 0
 is smaller or equal to $u_{+,a}$.
 The existence of such a solution will follow from the proof of Proposition \ref{voveru}.  In fact, from the proof
 of Proposition \ref{uv+-}, one can infer that there is a unique bounded solution to the equation
 $\frac D2 u''-r(x)u=0, \ -\infty<x<a$, with
$u(a)=1$.)
Similarly, for $a<0$,  let $u_{-,a}$ denote the solution to
\begin{equation}\label{u-a}
\begin{aligned}
&\frac D2 u''-r(x)u=0, \ a<x<\infty;\\
&u(a)=1;\\
&0\le u\le 1 \ \text{and}\ u\ \text{maximal}.
\end{aligned}
\end{equation}
For $a>0$, let $v_{+,a}$ denote the solution to
\begin{equation}\label{v+a}
\begin{aligned}
&\frac D2 v''-r(x)v=-1,\ -\infty<x<a;\\
&v(a)=0;\\
&v\ge0\ \text{and}\ v\ \text{minimal}.
\end{aligned}
\end{equation}
(The condition  that $v$ be minimal means that any other nonnegative solution is greater or equal to $v$.)
Similarly, for $a<0$  let $v_{-,a}$ denote the solution to
\begin{equation}\label{v-a}
\begin{aligned}
&\frac D2 v''-r(x)v=-1,\ a<x<\infty;\\
&v(a)=0;\\
&v\ge0\ \text{and}\ v\ \text{minimal}.
\end{aligned}
\end{equation}
\bf\noindent Remark.\rm\  It follows from the maximum principle that $u_{+,a}(x),u_{-,a}(x),v_{+,a}(x),v_{-,a}(x)$ are decreasing  in their dependence on  $r$.
Furthermore, if
 $r$ is sufficiently small, then $v_{+,a}$ on $(-\infty,a)$ and $v_{-,a}$ on $(a,\infty)$ will be equal to infinity.
The proof of Theorem \ref{comparison}-i shows that this occurs if
$r(x)\le\frac D{\gamma+x^2}$ for some $\gamma>0$ and sufficiently large $|x|$.

\begin{proposition}\label{voveru}
Let $r\gneqq0$ be a continuous function on $\mathbb{R}$. Then
for $a>0$,
$$
E_0^{(r)}T_a=\begin{cases}\infty,\ \text{if}\ v_{+,a}\equiv\infty;\\
\frac{v_{+,a}(0)}{u_{+,a}(0)},\ \text{otherwise},\end{cases}
$$
and for $a<0$,
$$
E_0^{(r)}T_a=
\begin{cases}\infty,\ \text{if}\ v_{-,a}\equiv\infty;\\
\frac{v_{-,a}(0)}{u_{-,a}(0)}, \ \text{otherwise}.\end{cases}
$$
\end{proposition}
\begin{proposition}\label{uv+-}
\noindent i. Let $\phi_3$ be as in Proposition \ref{criticality}. Then
\begin{equation}\label{u+3}
u_{+,a}(x)=\frac{\phi_3(x)}{\phi_3(a)}\frac{\int_{-\infty}^x\phi_3^{-2}(y)dy}{\int_{-\infty}^a\phi_3^{-2}(y)dy},\ x\le a;
\end{equation}
\begin{equation}\label{u-3}
u_{-,a}(x)=\frac{\phi_3(x)}{\phi_3(a)}\frac{\int_x^\infty\phi_3^{-2}(y)dy}{\int_a^\infty\phi_3^{-2}(y)dy},\ x\ge a;
\end{equation}
\begin{equation}\label{v+3}
v_{+,a}(x)=2\phi_3(x)\frac{\int_{-\infty}^xdy\thinspace\phi^{-2}_3(y)\int_x^adt\phi^{-2}_3(t)\int_y^t\phi_3(z)dz}
{\int_{-\infty}^a\phi^{-2}_3(y)dy},\ x\le a;
\end{equation}
\begin{equation}\label{v-3}
v_{-,a}(x)=2\phi_3(x)\frac{\int_x^\infty dy\thinspace\phi^{-2}_3(y)\int_a^xdt\phi^{-2}_3(t)\int_t^y\phi_3(z)dz}
{\int_a^\infty\phi^{-2}_3(y)dy},\ x\ge a.
\end{equation}
\noindent ii. Let $\phi_1$ be as in Proposition \ref{criticality}. Then
\begin{equation}\label{u+1}
u_{+,a}(x)=\frac{\phi_1(x)}{\phi_1(a)},\ x\le a;
\end{equation}
\begin{equation}\label{u-1}
u_{-,a}(x)=\frac{\phi_1(x)}{\phi_1(a)}\frac{\int_x^\infty\phi_1^{-2}(y)dy}{\int_a^\infty\phi_1^{-2}(y)dy},\ x\ge a;
\end{equation}
\begin{equation}\label{v+1}
v_{+,a}(x)=2\phi_1(x)\int_x^ady\thinspace\phi^{-2}_1(y)\int_{-\infty}^y\phi_1(z)dz,\ x\le a;
\end{equation}
\begin{equation}\label{v-1}
v_{-,a}(x)=2\phi_1(x)\frac{\int_x^\infty dy\thinspace\phi^{-2}_1(y)\int_a^xdt\phi^{-2}_1(t)\int_t^y\phi_1(z)dz}
{\int_a^\infty\phi^{-2}_1(y)dy},\ x\ge a.
\end{equation}
\end{proposition}
\bf\noindent Remark 1.\rm\ Of course, formulas similar to \eqref{u+1}-\eqref{v-1} can be given in terms of $\phi_2$.

\bf\noindent Remark 2.\rm\  As noted in the remark after \eqref{v-a}, when $r$ is sufficiently small, $v_{+,a}$ and $v_{-,a}$ are infinite. For $v_{+,a}$, one sees from
\eqref{v+3} that the infiniteness is equivalent to $\int_{-\infty}\phi^{-2}_3(y)\int_y^0\phi_3(z)dz=\infty$, and from \eqref{v+1} that it is equivalent
to $\int_{-\infty}\phi_1(y)dy=\infty$. Similar equivalences hold for the infiniteness of $v_{-,a}$ from \eqref{v-3} and \eqref{v-1}.

As an immediate corollary of Propositions \ref{voveru} and Proposition \ref{uv+-}
we obtain the following explicit representation of $E_0^{(r)}T_a$.
\begin{theorem}\label{rep}
\noindent i. Let $\phi_3$ be as in Proposition \ref{criticality}.
Then
\begin{equation}\label{withphi3}
E_0^{(r)}T_a=\begin{cases}\frac{2\phi_3(a)}{\int_{-\infty}^0\phi_3^{-2}(x)dx}
\int_{-\infty}^0dx\thinspace\phi_3^{-2}(x)\int_0^ady\thinspace\phi_3^{-2}(y)\int_x^y\thinspace\phi_3(z)dz,\ a>0;\\
\frac{2\phi_3(a)}{\int_0^{\infty}\phi_3^{-2}(x)dx}\int_0^{\infty}dx\thinspace\phi_3^{-2}(x)
\int_a^0dy\thinspace\phi_3^{-2}(y)\int_y^x\thinspace\phi_3(z)dz,\ a<0.
\end{cases}
\end{equation}
\noindent ii. Let $\phi_1$ be as in Proposition \ref{criticality}.
Then
\begin{equation}\label{withphi1}
E_0^{(r)}T_a=\begin{cases}2\phi_1(a)\int_0^ady\thinspace\phi_1^{-2}(y)\int_{-\infty}^y\phi_1(z)dz,\ a>0;\\
\frac{2\phi_1(a)}{\int_0^{\infty}\phi_1^{-2}(x)dx}\int_0^{\infty}dx\thinspace\phi_1^{-2}(x)
\int_a^0dy\thinspace\phi_1^{-2}(y)\int_y^x\thinspace\phi_1(z)dz,\ a<0.
\end{cases}
\end{equation}
\end{theorem}
\bf\noindent Remark 1.\rm\ In light of Remark 1 after Proposition \ref{uv+-}, a formula analogous to \eqref{withphi1} holds with $\phi_2$ in  place of $\phi_1$.

\bf\noindent Remark 2.\rm\ In the case that $r(x)=r>0$ is constant, letting $\phi_3(x)=\exp(\sqrt{\frac2Dr}\thinspace x)+\exp(-\sqrt{\frac2Dr}\thinspace x)$
and $\phi_1(x)=\exp(\sqrt{\frac2Dr}\thinspace x)$, one can check that
\eqref{withphi3} and \eqref{withphi1} yield \eqref{constantr}.

Using Theorem \ref{rep} with Propositions \ref{voveru} and \ref{uv+-}, along with the fact that  the functions $u_{+,a},u_{-,a}$, $v_{+,a},v_{-,a}$  are decreasing in
$r$, and choosing  test functions $\phi_3$ appropriately, we will  prove the following quantitative estimates on  $E_0^{(r)}T_a$ in terms of the behavior of the resetting rate $r$.
\begin{theorem}\label{comparison1}
Let $l>-1$.
If
\begin{equation}\label{rgrowth}
c_1(\gamma_1+x^2)^l\le r(x)\le c_2(\gamma_2+x^2)^l, \ \text{for all}\ x\in\mathbb{R},  \text{where}\ c_2\ge c_1>0\ \text{and}\ \gamma_2,\gamma_1>0,
\end{equation}
then there exist $K_2>K_1>0$ and
$M_1,M_2>0$ such that
\begin{equation}\label{solugrowth}
|a|M_1e^{K_1(1+a^2)^{\frac{l+1}2}}\le E_0^{(r)}T_a\le M_2e^{K_2(1+a^2)^{\frac{l+1}2}},\ \text{for all}\ a\in\mathbb{R}.
\end{equation}
\end{theorem}
\bf\noindent Remark.\rm\ It follows that if $c_1(\gamma_1+x^2)^l\le r(x)\le c_2(\gamma_2+x^2)^l$, for some $l>-1$, then a necessary condition for the finiteness of
$\int_{\mathbb{R}}E_0^{(r)}T_a\mu(da)$ is that all  the moments of $\mu$ are finite.

\begin{theorem}\label{comparison}
\noindent i. If $r(x)\le \frac D{\gamma+x^2}$, for some $\gamma>0$ and for sufficiently large $|x|$, then $E_0^{(r)}T_a=\infty$, for all $a\in\mathbb{R}-\{0\}$;

\noindent ii. If $r(x)\ge\frac{D\lambda}{\gamma+x^2}$, for some $\lambda>1$ and some $\gamma>0$, and for sufficiently large $|x|$, then
$E_0^{(r)}T_a<\infty$, for all $a\in\mathbb{R}$.

\noindent iii.
For any $\lambda>1$, there exists an  $r(x)$ satisfying
$r(x)\sim \frac{D\lambda}{x^2}$ as $x\to\infty$, and such that  $E_0^{(r)}T_a\sim C|a|^{\frac{1+\sqrt{1+8\lambda}}2}$,  as $|a|\to\infty$, for some $C>0$.

\noindent iv. If
$$
\frac{D\lambda_1}{\gamma_1+x^2}\le r(x)\le \frac{D\lambda_2}{\gamma_2+x^2}, \ \text{for all}\ x, \ \text{where}\
1<\lambda_1\le\lambda_2\ \text{and}\ \gamma_1,\gamma_2>0,
$$
then for any $\epsilon>0$, there exist $C_1,C_2>0$ such that
\begin{equation}\label{witheps}
C_1|a|^{\frac{1+\sqrt{1+8\lambda_1}}2-\epsilon }\le E_0^{(r)}T_a\le C_2|a|^{\frac{1+\sqrt{1+8\lambda_2}}2+\epsilon},\ \text{for}\ |a|\ge1.
\end{equation}
\end{theorem}
\noindent \bf Remark.\rm\ We expect that part (iv) also holds with $\epsilon=0$.


\medskip

Theorems \ref{comparison1} and \ref{comparison}  show generically that if $r(x)$ is of the order $|x|^{2l}$, with $l>-1$, then
$\log E_0^{(r)}T_a$ is of the order $|a|^{l+1}$; in particular, the smaller the asymptotic size of $r$, the smaller the asymptotic growth
rate of $E_0^{(r)}T_a$. The asymptotic growth rate of $E_0^{(r)}T_a$ continues to decrease when $r(x)\sim \frac{D\lambda}{x^2}$
with $\lambda>1$; now the growth rate of $E_0^{(r)}T_a$ is more or less of the order $|a|^{\frac{1+\sqrt{1+8\lambda}}2}$. Note that this exponent
increases to $\infty$ when $\lambda$ increases to $\infty$ and decreases to 2 when $\lambda$ decreases to 1. However, if $\lambda=1$, then
 $E_0^{(r)}T_a=\infty$, for all $a\neq0$.

Theorem \ref{comparison} shows that the dependence of
$E_0^{(r)}T_a$ on $r$ is very sensitive when $r$ has quadratic decay.
If one uses regularly varying resetting rates $r$,  Theorem \ref{comparison} shows that
 if $\mu$ is such that its  $p$th moment ($p>0$) is finite if and only if $p<p_0$,
then if $p_0\le 2$, $\int_{\mathbb{R}}(E_0^{(r)}T_a)\thinspace\mu(da)$ will be infinite for all  $r$,
and if $p_0>2$, then $\int_{\mathbb{R}}(E_0^{(r)}T_a)\thinspace\mu(da)$  will be finite
only for a very narrow window of these rates;
namely, for resetting rates $r(x)$ that satisfy $r(x)\sim \frac{D\lambda}{x^2}$, where
$1<\lambda<\frac{(2p_0-1)^2-1}8$.

Theorems \ref{comparison1} and \ref{comparison} would seem to  suggest that for many distributions $\mu$ supported on all of $\mathbb{R}$, a near-optimal (or optimal) $r$ for which
$\int_{\mathbb{R}}(E_0^{(r)}T_a)\thinspace\mu(da)$ will be close to minimal (or minimal) will be one with quadratic decay.
We have some numerical work that doesn't quite bear this out, however it doesn't disprove this thesis.
We consider  the case that the target distribution $\mu$ is a two-sided symmetric exponential distribution: $\mu((x,\infty))=\frac12e^{-\beta x}$, for $x>0$ and some $\beta>0$, and $\mu((-\infty,-x))=\mu((x,\infty))$, for $x\ge0$.
We compare
the expected time to locate the target in the case that the optimal
\it constant\rm\ resetting rate is used to the case that  certain \it quadratically decaying\rm\ resetting rates are used.
By \eqref{constantr}, the expected time to locate the target with constant resetting rate $r>0$ is
$$
\int_0^\infty\frac{e^{\sqrt{\frac{2r}D}\thinspace a}-1}r\beta e^{-\beta a}da.
$$
A standard calculation reveals that this expression is minimized when $r=\frac{D\beta^2}8$, and that the minimum value
is $\frac8{D\beta^2}$. Recall that the expected distance to the target is $\frac1\beta$; thus the optimal
expected time to locate the target is $\frac8D$ times the square of the expected distance to the target.

Now consider $r(x)=\frac{m(m-1)}2D\frac{|x|^{m-2}}{\gamma+|x|^m}$, where $m>2$ and $\gamma>0$.
(For any $\gamma>0$, and with $\lambda=\frac{m(m-1)}2$, this function is an example of the function $r$
appearing in Theorem \ref{comparison}-iii.)
The function $\phi_3(x)=\gamma+|x|^m$ satisfies $\frac D2\phi_3''(x)-r(x)\phi_3(x)=0$.
Thus $E_0^{(r)}T_a$ is given by \eqref{withphi3} with this choice of $\phi_3$.
Since $r$ is symmetric,  $E_0^{(r)}T_a$ is symmetric in $a$, and thus the expected time to locate the target is
$\int_0^\infty (E_0^{(r)}T_a)\beta e^{-\beta a}da$.
We substitute   in this integral the expression for   $E_0^{(r)}T_a$ in \eqref{withphi3}.
We want to minimize the resulting quantity as $\gamma$ varies over $(0,\infty)$ and $m$ varies over $(2,\infty)$.
My colleague, Nir Gavish, found that
for $\beta\ge0.04$, the infimum value as a function of $\beta$ can be approximated by the function $\frac{8.14+12.42\exp(-35.66\beta)}{D\beta^2}$, with an error of less than 1 percent in the numerator.
This is slightly worse than what we obtained using the optimal constant resetting rate.
Of course, it is still possible that the infimum of $\int_0^\infty (E_0^{(r)}T_a)\beta e^{-\beta a}da$
over all $r$ of the form $r(x)=\frac{c_1}{c_2+x^2}$, with $c_1,c_2>0$  is less than
$\frac8{D\beta^2}$. And we certainly expect that the infimum over \it all\rm\ $r$ that exhibit quadratic decay will be less
than $\frac8{D\beta^2}$.

The discussion above, as well as our results, have been geared in particular to the case that the support of the target
distribution is all of $\mathbb{R}$. If the support of the target distribution is, say,
$[-L_1,L_2]$, where $L_1,L_2>0$,
 then there is no reason to search outside of this interval, and thus as soon as the searcher reaches
$-L_1$ or $L_2$, its position should be reset to 0. This is equivalent to setting $r\equiv\infty$ off of $[-L_1,L_2]$.
We discuss this situation  in section \ref{finsupp}.

We end this presentation  of results by noting that Proposition \ref{criticality} and Theorem \ref{rep} furnish
explicit, albeit rather complicated, variational formulas for $\inf_{r\gneqq0}\int_{\mathbb{R}}(E_0^{(r)}T_a)\mu(da)$.
Assume that $\mu$ has mass both in $(0,\infty)$ and in $(-\infty,0)$, and for  convenience,  assume that the origin is not an atom of the distribution $\mu$. Then  $\mu$ can be written in the form
\begin{equation*}
\begin{aligned}
&\mu=(1-p)\mu_-+\thinspace p\mu_+,\
\text{where}\ p\in(0,1), \
 \text{and}\ \mu_+\ \text{and}\ \mu_- \\
 & \text{are  probability measures on}\ (0,\infty)\ \text{and}\ (-\infty, 0)\ \text{respectively}.
\end{aligned}
\end{equation*}
\begin{corollary}\label{variational}
\noindent i.
\begin{equation}
\begin{aligned}
&\inf_{0\lneqq r\in C(\mathbb{R})}\int_{-\infty}^\infty E^{(r)}_0T_a\thinspace\mu(da)=
\inf_{\stackrel{\phi\in C^2(\mathbb{R}),\phi>0,\phi''\gneqq0}{\int_{-\infty}^\infty\phi^{-2}(x)dx<\infty}}\\
&\Big[\frac {2p}{\int_{-\infty}^0\phi^{-2}(x)dx}\int_0^\infty \mu_+(da)\phi(a)\int_{-\infty}^0dx\thinspace\phi^{-2}(x)\int_0^ady\thinspace\phi^{-2}(y)\int_x^ydz\thinspace\phi(z)+\\
&\frac{2(1-p)}{\int_0^{\infty}\phi^{-2}(x)dx}\int_{-\infty}^0\mu_-(da)\phi(a)\int_0^{\infty}dx\thinspace\phi^{-2}(x)
\int_a^0dy\thinspace\phi^{-2}(y)\int_y^xdz\thinspace\phi_3(z)\Big].
\end{aligned}
\end{equation}
\noindent ii.
\begin{equation}
\begin{aligned}
&\inf_{0\lneqq r\in C(\mathbb{R})}\int_{-\infty}^\infty E^{(r)}_0T_a\thinspace\mu(da)=
\inf_{\stackrel{\phi\in C^2(\mathbb{R}),\phi>0,\phi''\gneqq0}{\int_{-\infty}\phi^{-2}(x)dx=\infty}}\\
&\Big[2p\int_0^\infty \mu_+(da)\phi(a)\int_0^ady\thinspace\phi^{-2}(y)\int_{-\infty}^y\phi(z)dz+\\
&\frac{2(1-p)}{\int_0^{\infty}\phi^{-2}(x)dx}\int_{-\infty}^0\mu_-(da)\phi(a)\int_0^{\infty}dx\thinspace\phi^{-2}(x)
\int_a^0dy\thinspace\phi^{-2}(y)\int_y^xdz\thinspace\phi(z)\Big].
\end{aligned}
\end{equation}
\end{corollary}
\bf \noindent Remark.\rm\ If $\phi>0,\phi''\gneqq0$ and $\int_{-\infty}\phi^{-2}(x)dx=\infty$, then necessarily
$\int^{\infty}\phi^{-2}(x)dx<\infty$ (see the proof of Proposition \ref{criticality}), so there is no need to include this last condition in part (ii).

Consider the case that $\mu$ is symmetric; that is, the case that $p=\frac12$ and $\mu_+(A)=\mu_-(-A)$, for $A\subset(0,\infty)$.
Then presumably,
$$
\inf_{0\lneqq r\in C(\mathbb{R})}\int_{-\infty}^\infty E^{(r)}_0T_a\thinspace\mu(da)=
\inf_{\stackrel{0\lneqq r\in C(\mathbb{R})}{r \ \text{is even}}}\int_{-\infty}^\infty E^{(r)}_0T_a\thinspace\mu(da),
$$
although we don't have a proof.
\begin{corollary}\label{variationalsym}
Assume that $\mu$ is symmetric.

\noindent i.
\begin{equation}
\begin{aligned}
&\inf_{\stackrel{0\lneqq r\in C(\mathbb{R})}{r\ \text{\rm is even}}}\int_{-\infty}^\infty E^{(r)}_0T_a\thinspace\mu(da)=
\inf_{\stackrel{\phi\in C^2(\mathbb{R}),\phi>0,\phi''\gneqq0, \phi\ \text{\rm is even}}{\int_{-\infty}^\infty\phi^{-2}(x)dx<\infty}}\\
&\Big[\frac {2}{\int_{-\infty}^0\phi^{-2}(x)dx}\int_0^\infty \mu_+(da)\phi(a)\int_{-\infty}^0dx\thinspace\phi^{-2}(x)\int_0^ady\thinspace\phi^{-2}(y)\int_x^ydz\thinspace\phi(z)].
\end{aligned}
\end{equation}
\noindent ii.
\begin{equation}
\begin{aligned}
&\inf_{\stackrel{0\lneqq r\in C(\mathbb{R})}{r\ \text{\rm is even}}}\int_{-\infty}^\infty E^{(r)}_0T_a\thinspace\mu(da)=
\inf_{\stackrel{\phi\in C^2(\mathbb{R},\phi>0,\phi''\gneqq0, \frac{\phi''}\phi \text{\rm\ is even}}{\int_{-\infty}\phi^{-2}(x)dx=\infty}}\\
&\Big[2\int_0^\infty \mu_+(da)\phi(a)\int_0^ady\thinspace\phi^{-2}(y)\int_{-\infty}^y\phi(z)dz\Big].
\end{aligned}
\end{equation}
\end{corollary}
\bf\noindent Remark.\rm\ In part (ii), $\phi$ cannot be even because, as noted in the remark following Corollary \ref{variational},
the conditions $\phi>0,\phi''\gneqq0$ and $\int_{-\infty}\phi^{-2}(x)dx=\infty$ lead automatically to
$\int^{\infty}\phi^{-2}(x)dx<\infty$.

We prove Propositions \ref{criticality}-\ref{uv+-} in sections \ref{critical}-\ref{proofuv+-} respectively, and Theorems \ref{comparison1} and \ref{comparison} in sections \ref{proofthmcomp1} and \ref{proofthmcomp} respectively.
In section \ref{finsupp} we discuss the case in which the target distribution is supported on a finite interval.

\section{Proof of Proposition \ref{criticality}}\label{critical}
The proof is an application of the criticality theory of second order elliptic operators---see \cite[chapter 4]{P}.
The operator $L:=\frac D2\frac{d^2}{dx^2}-r(x)$ on $\mathbb{R}$ with $r\gneqq0$ is subcritical, and thus there exists a positive $L$-harmonic function $\phi$ (that is, $\phi>0$ and $L\phi=0$).
Denote by $L^\phi$  the operator which is the  $h$-transform of $L$ by the function $\phi$;  $L^\phi u:=\frac1\phi L(\phi u)$.
It follows that $L^\phi=\frac D2\frac{d^2}{dx^2}+D\frac{\phi'}\phi\frac d{dx}$.
Subcriticality is preserved by $h$-transforms, so $L^\phi$ is also subcritical.  For one-dimensional operators, in the subcritical case, the cone of positive
harmonic functions is two-dimensional. Furthermore, one of them is minimal at $-\infty$ and the other one is minimal at $+\infty$.
Denote by $\hat\phi_1$ and $\hat\phi_2$ positive $L^\phi$-harmonic functions, with $\hat\phi_1$ minimal at $-\infty$ and $\hat\phi_2$ minimal at $+\infty$.
Define $\phi_1=\phi\hat\phi_1$ and $\phi_2=\phi\hat\phi_2$. Then $\phi_1$ and $\phi_2$ are $L$-harmonic; indeed,
$0=(L^\phi)(\hat\phi_i)=\frac1{\phi}L(\phi\hat\phi_i)$, $i=1,2$.
Since the zeroth order term in $L^\phi$ vanishes, $L^\phi$ is the generator of a diffusion process. For diffusion generators, subcriticality is equivalent to transience.
By the Martin boundary theory (\cite[chapter 7]{P}, the $h$-transformed operators $(L^\phi)^{\hat\phi_1}$ and  $(L^\phi)^{\hat\phi_2}$ are transient to $+\infty$ and $-\infty$
respectively. These operators are given by
$$
(L^\phi)^{\hat\phi_1}=\frac D2\frac{d^2}{dx^2}+D\frac{\phi'}\phi \frac d{dx}+D\frac{\hat\phi'_1}{\hat\phi_1}\frac d{dx}
$$
and
$$
(L^\phi)^{\hat\phi_2}=\frac D2\frac{d^2}{dx^2}+D\frac{\phi'}\phi \frac d{dx}+D\frac{\hat\phi'_2}{\hat\phi_2}\frac d{dx}.
$$
The transience to $+\infty$ ($-\infty$) of a diffusion generator $\frac D2\frac{d^2}{dx^2}+Db(x)\frac d{dx}$ is equivalent
to $\int_{-\infty}\exp(-\int_0^x2b(y)dy)dx=\infty$ and $\int^{\infty}\exp(-\int_0^x2b(y)dy)dx<\infty$
($\int_{-\infty}\exp(-\int_0^x2b(y)dy)dx<\infty$ and $\int^{\infty}\exp(-\int_0^x2b(y)dy)dx=\infty$).
Since $b=\frac{\phi'}\phi+\frac{\hat\phi_i'}{\hat\phi_i}$ for $(L^\phi)^{\hat\phi_i}$, $i=1,2$,
 it follows that
 $\phi_1$ and $\phi_2$ satisfy \eqref{intcond}. Define $\phi_3=\phi_1+\phi_2$. Then $\phi_3$ is also $L$-harmonic and it satisfies \eqref{intcond}. If $r$ is an even function, then $\phi_3(-x)$ is also $L$-harmonic. Thus,
$\bar\phi_3(x):=\phi_3(x)+\phi_3(-x)$ is $L$-harmonic, satisfies \eqref{intcond} (for $\phi_3$) and is even.
\hfill $\square$

\section{Proof of Proposition \ref{voveru}}\label{proofpropvoveru}
We will prove the proposition for $a>0$; the same type of proof works for $a<0$.
For $n>0$, let $\mathcal{T}_t^{(n)}$ be the semigroup defined by
$$
\mathcal{T}_t^{(n)}f(x)=E_x^{(r)}(f(X(t));T_a\wedge T_{-n}>t):=E_x^{(r)}(f(X(t))1_{\{T_a\wedge T_{-n}>t\}}),\ x\in[-n,a],
$$
for bounded continuous $f$.
Its generator is $\mathcal{L}$ as in \eqref{generator} with the zero Dirichlet boundary condition at $x=a$ and at $x=-n$.
Let $w_n(x,t)=\mathcal{T}_t^{(n)}1(x)$. Then $w_n(x,t)=P_x^{(r)}(T_a\wedge T_{-n}>t)$ and it solves
\begin{equation}\label{w}
\begin{aligned}
&\frac{\partial}{\partial t}w_n=\mathcal{L}w_n=\frac D2w_n''+r(x)\big(w_n(0,t)-w_n(x,t)\big), \ x\in(-n,a);\\
&w_n(x,0)=1,\ x\in(-n,a);\ \
w_n(a,t)=w_n(-n,t)=0,\ t>0.
\end{aligned}
\end{equation}
Let
$$
A_n(x,s)=\int_0^\infty\exp(-st)w_n(x,t)dt,\ x\in[-n,a],\ s>0.
$$
Then
$$
\begin{aligned}
&\frac D2\frac{d^2A_n}{dx^2}(x,s)+r(x)\big(A_n(0,s)-A_n(x,s)\big)=\mathcal{L}A_n(x,s)=\\
&\int_0^\infty\exp(-st)\mathcal{L}w_n(x,t)dt=\int_0^\infty\exp(-st)\frac{\partial}{\partial t}w_n(x,t)dt=\\
&-1+s\int_0^\infty \exp(-st)w_n(x,t)dt=-1+sA_n(x,s),\ \text{for}\ x\in(-n,a).
\end{aligned}
$$
Letting $s\to0$, we find that $A_n(x):=A_n(x,0)$ satisfies
\begin{equation}\label{A}
\begin{aligned}
&\frac D2A_n''(x)-r(x)A_n(x)=-1-r(x)A_n(0),\ x\in(-n,a);\\
&A_n(a)=A_n(-n)=0.
\end{aligned}
\end{equation}
Note that
$$
A_n(x)=\int_0^\infty w_n(x,t)dt=\int_0^\infty P_x^{(r)}(T_a\wedge T_{-n}>t)dt=E_x^{(r)}T_a\wedge T_{-n};
$$
thus
\begin{equation}\label{A0}
E_0^{(r)}T_a\wedge T_{-n}=A_n(0).
\end{equation}

For $c>0$, let
$B_{n,c}(x)$ solve the equation
\begin{equation}
\begin{aligned}
&\frac D2B_{n,c}''-r(x)B_{n,c}(x)=-1-cr(x),\ x<a;\\
&B_{n,c}(a)=B_{n,c}(-n)=0.
\end{aligned}
\end{equation}
We look for $c_n>0$ satisfying $B_{n,c_n}(0)=c_n$.
It then follows that $A_n(x)=B_{n,c_n}(x)$, and in particular,
\begin{equation}\label{A0c0}
A_n(0)=c_n.
\end{equation}
Let $v_{n,+,a}$ solve  the equation
\begin{equation}\label{vn+a}
\begin{aligned}
&\frac D2 v''_{n,+,a}-r(x)v_{n,+,a}=-1,\ x\in(-n,a);\\
&v_{n,+,a}(a)=v_{n,+,a}(-n)=0,
\end{aligned}
\end{equation}
and let $u_{n,+,a}$ solve the equation
\begin{equation}\label{un+a}
\begin{aligned}
&\frac D2u_{n,+,a}''-r(x)u_{n,+,a}=0,\ x\in(-n,a);\\
&u_{n,+,a}(a)=u_{n,+,a}(-n)=1.
\end{aligned}
\end{equation}
(We note that by the maximum principle, $0\le u_{n,+,a}\le 1$.)
Then $B_{n,c}=v_{+,a}+c(1-u_{n,+,a})$, and thus the equation  $B_{n,c_n}(n)=c_n$ is solved by
$c_n=\frac{v_{n,+,a}(0)}{u_{n,+,a}(0)}$.
 Thus, from \eqref{A0} and \eqref{A0c0},
\begin{equation}\label{finaln}
 E_0^{(r)}T_a\wedge T_{-n}=A_n(0)=c_n=\frac{v_{n,+,a}(0)}{u_{n,+,a}(0)}.
\end{equation}
By the maximum principle $\lim_{n\to\infty}u_{n,+,a}=u_{+,a}$ and
$\lim_{n\to\infty}v_{n,+,a}=v_{+,a}$, where $u_{+,a}$ and $v_{+,a}$ are given by \eqref{u+a} and \eqref{v+a}. Thus, letting $n\to\infty$ in \eqref{finaln} gives
$E_0^{(r)}T_a=\infty$, if $v_{+,a}\equiv\infty$;  otherwise it gives
$$
E_0^{(r)}T_a=\frac{v_{+,a}(0)}{u_{+,a}(0)}.
$$
\hfill $\square$

\section{Proof of Proposition \ref{uv+-}}\label{proofuv+-}
  We'll prove the formulas for $u_{+,a}$ and $v_{+,a}$. The proofs for $u_{-,a}$ and $v_{-,a}$ are similar.
 Let $L:=\frac D2\frac{d^2}{dx^2}-r(x)$.

The solution $u_{+,a}$ to \eqref{u+a} is obtained as
$$
u_{+,a}=\lim_{N\to\infty}u_{+,a,N},
$$
where $u_{a,+,N}$ satisfies $Lu_{a,+,N}=0$ in $(-N,a)$ with boundary condition $u_{+,a,N}(a)=u_{+,a,N}(-N)=1$.
Let $\phi_i$ be as in Proposition \ref{criticality}, with either $i=1$ or $i=3$. As noted in section \ref{critical}, the $h$-transform of $L$ by $\phi_i$, denoted by $L^{\phi_i}$, is defined by $L^{\phi_i} u=\frac1{\phi_i} L(\phi_i u)$,
and when written out, one obtains $L^{\phi_i}=\frac D2\frac{d^2}{dx^2}+D\frac{\phi_i'}{\phi_i}\frac d{dx}$.
Write $u_{+,a,N}$  in the form $u_{+,a,N}=\phi_i\bar u_{+,a,N}$.
Since $Lu_{+,a,N}=0$, we have
$$
0=L^{\phi_i}\bar u_{+,a,N}=\frac D2\bar u_{+,a,N}''+D\frac{\phi_i'}{\phi_i}\bar u_{+,a,N}'=\frac D2\frac1{\phi_i^2}\big(\phi_i^2\overline{u}_{+,a,N}'\big)',
$$
and since $u_{+,a,N}(a)=u_{+,a,N}(-N)=1$,
we have $\bar u_{+,a,N}(a)=\frac1{\phi_i(a)}$ and $\bar u_{+,a,N}(-N)=\frac1{\phi_i(-N)}$.
Solving by integrating twice and using the boundary condition, we obtain
$$
\bar u_{+,a,N}(x)=\frac1{\phi_i(a)}-\big(\frac1{\phi_i(a)}-\frac1{\phi_i(-N)}\big)
\frac{\int_x^a\phi_i^{-2}(y)dy}{\int_{-N}^a\phi_i^{-2}(y)dy}.
$$

Choose first $i=3$. Since $\int_{-\infty}\phi_i^{-2}(y)dy<\infty$, it follows that there exists a sequence
$N_k\stackrel{k\to\infty}{\to}\infty$ such that $\lim_{k\to\infty}\phi_3(-N_k)=\infty$.
Thus substituting $N_k$ for $N$ above and letting $k\to\infty$,
 we conclude that $u_{+,a}$ is given by \eqref{u+3}.
 (Retroactively, it then follows from the uniqueness of the solution to \eqref{u+a} that in fact $\lim_{N\to\infty}\phi_3(-N)=\infty$.)
Now choose $i=1$.
Letting
$N\to\infty$, we conclude that $u_{+,a}$ is given  by
\eqref{u+1}.

The solution $v_{+,a}$ to \eqref{v+a} is obtained as
$$
\lim_{N\to\infty}v_{+,a,N},
$$
where $v_{+,a,N}$ solves $Lv_{+,a,N}=-1$ in $(-N,a)$ with boundary condition $v_{+,a,N}(a)=v_{+,a,N}(-N)=0$.
As was done above in solving for $u_{+,a}$, we make an $h$-transform with $\phi_i$, $i=1,3$.
Similar to the above, we write $v_{+,a,N}$ in the form $v_{+,a,N}=\phi_i\bar v_{+,a,N}$, and then obtain the
equation
$$
-\frac1{\phi_i}=\frac D2\bar v_{+,a,N}+D\frac{\phi_i'}{\phi_i}\bar v_{+,a,N}'=\frac D2\frac1{\phi_i^2}\big(\phi_i^2\overline{v}_{+,a,N}'\big)',
$$
with the boundary condition $\bar v_{+,a,N}(a)=\bar v_{+,a,N}(-N)=0$.
Solving by integrating twice and using the boundary condition, we obtain
$$
\bar v_{+,a,N}(x)=c_N\int_x^a\phi_i^{-2}(y)dy-2\int_x^a dy\phi_i^{-2}(y)\int_y^a\phi_i(z)dz,
$$
where
$$
c_N=\frac{2\int_{-N}^a dy\phi_i^{-2}(y)\int_y^a\phi_i(z)dz}{\int_{-N}^ady\phi_i^{-2}(y)}.
$$

Choosing $i=1$, we have $\lim_{N\to\infty}c_N=2\int_{-\infty}^a\phi_1(y)dy$, and thus
$$
\begin{aligned}
&v_{+,a}(x)=\lim_{N\to\infty}v_{+,a,N}(x)=\phi_1(x)\lim_{N\to\infty}\bar v_{+,a,N}(x)=\\
&2\phi_1(x)\int_x^ady\thinspace\phi^{-2}_1(y)\int_{-\infty}^y\phi_1(z)dz,
\end{aligned}
$$
as in \eqref{v+1}.

Choosing $i=3$, we have
\begin{equation}\label{squarebracket}
\begin{aligned}
&v_{+,a}(x)=\lim_{N\to\infty}v_{+,a,N}(x)=\phi_3(x)\lim_{N\to\infty}\bar v_{+,a,N}(x)=\frac{2\phi_3(x)}{\int_{-\infty}^a\phi_3^{-2}(y)dy}\times\\
&\Big[\big(\int_{-\infty}^a dy\phi_3^{-2}(y)\int_y^a\phi_3(z)dz\big)
\big(\int_x^a\phi_3^{-2}(t)dt\big)-\\
&\big(\int_x^a dy\phi_3^{-2}(y)\int_y^a\phi_3(z)dz\big)
\big(\int_{-\infty}^a\phi_3^{-2}(t)dt\big)\Big].
\end{aligned}
\end{equation}
Write the first term of the square brackets in \eqref{squarebracket}   as
\begin{equation*}
\begin{aligned}
&\big(\int_{-\infty}^a dy\phi_3^{-2}(y)\int_y^a\phi_3(z)dz\big)
\big(\int_x^a\phi_3^{-2}(t)dt\big)=\\
&\big(\int_{-\infty}^x dy\phi_3^{-2}(y)\int_y^a\phi_3(z)dz\big)
\big(\int_x^a\phi_3^{-2}(t)dt\big)+\\
&\big(\int_x^a dy\phi_3^{-2}(y)\int_y^a\phi_3(z)dz\big)
\big(\int_x^a\phi_3^{-2}(t)dt\big),
\end{aligned}
\end{equation*}
and write the second term there as
\begin{equation}\label{secondterm}
\begin{aligned}
&\big(\int_x^a dy\phi_3^{-2}(y)\int_y^a\phi_3(z)dz\big)
\big(\int_{-\infty}^a\phi_3^{-2}(t)dt\big)=\\
&\big(\int_x^a dy\phi_3^{-2}(y)\int_y^a\phi_3(z)dz\big)
\big(\int_{-\infty}^x\phi_3^{-2}(t)dt\big)+\\
&\big(\int_x^a dy\phi_3^{-2}(y)\int_y^a\phi_3(z)dz\big)
\big(\int_x^a\phi_3^{-2}(t)dt\big).
\end{aligned}
\end{equation}
Notice that the final line from each of  the preceding two displays is the same.
Thus, the expression in the square brackets in \eqref{squarebracket} is equal to
\begin{equation}\label{squarebracketagain}
\begin{aligned}
&\big(\int_{-\infty}^x dy\phi_3^{-2}(y)\int_y^a\phi_3(z)dz\big)
\big(\int_x^a\phi_3^{-2}(t)dt\big)-\\
&\big(\int_x^a dt\phi_3^{-2}(t)\int_t^a\phi_3(z)dz\big)
\big(\int_{-\infty}^x\phi_3^{-2}(y)dy\big).
\end{aligned}
\end{equation}
(The second term above is the middle line of \eqref{secondterm}, but we have switched the roles of the variables of integration $t$ and $y$.)
The expression in \eqref{squarebracketagain} can be written as
$$
\int_{-\infty}^xdy\thinspace\phi^{-2}_3(y)\int_x^adt\phi^{-2}_3(t)\int_y^t\phi_3(z)dz.
$$
Substituting this for the expression in the square brackets in \eqref{squarebracket}, we conclude that $v_{+,a}$ is given by \eqref{v+3}.
\hfill $\square$

\section{Proof of Theorem \ref{comparison1}}\label{proofthmcomp1}
  We will prove the estimate for $a>0$; the same type of proof
works for $a<0$. Fix $l>-1$ and assume that $r$ satisfies \eqref{rgrowth}.
That is, let
$$
r_-(x)=c_1(\gamma_1+x^2)^l,\ \ \  \ r_+(x)=c_2(\gamma_2+x^2)^l,
$$
where $0<c_1<c_2$ and $\gamma_1,\gamma_2>0$, and where   $r_-(x)\le r(x)\le r_+(x)$, for all $x\in\mathbb{R}$.
We will need to compare the solutions $u_{+,a}$ and $v_{+,a}$ of \eqref{u+a}  and \eqref{v+a}
 for different choices of the function $r$, so we will denote them here by $u_{+,a,r}$ and $v_{+,a,r}$.
By Proposition \ref{voveru} and the fact that $u_{+,a,r}(0)$ and $v_{+,a,r}(0)$ are  decreasing
in their dependence on $r$, it follows that
\begin{equation}\label{rr+r-}
\frac{v_{+,a,r_+}(0)}{u_{+,a,r_-}(0)}\le E_0^{(r)}T_a\le\frac{v_{+,a,r_-}(0)}{u_{+,a,r_+}(0)}.
\end{equation}

Define
\begin{equation}\label{phis}
\phi_3(x)=e^{\psi(x)},\ \text{where}\ \psi(x)=\lambda(\gamma+x^2)^{\frac{l+1}2},\ \text{with}\ \gamma,\lambda>0.
\end{equation}
Let
\begin{equation}\label{hatr}
\begin{aligned}
&\hat r(x):=\frac D2\frac{\phi_3''(x)}{\phi_3(x)}=\frac D2\big((\psi'(x))^2+\psi''(x)\big)=\\
&\frac D2\lambda(l+1)(\gamma+x^2)^{\frac{l-3}2}
\big[(l+1)\lambda x^2(\gamma+x^2)^{\frac{l+1}2}+\gamma+lx^2\big].
\end{aligned}
\end{equation}
We will show that one can choose $\gamma$ and $\lambda$ so that the corresponding $\hat r$, which we will denote by
$\hat r_-$, satisfies
 $0\le \hat r_-\le r_-$, and we will show that one can choose
$\gamma$ and $\lambda$ so that the corresponding $\hat r$, which we will denote by
$\hat r_+$, satisfies
 $\hat r_+\ge r_+$.
\it  Note that since \eqref{hatr} does not depend on $a$, the $\gamma$ and $\lambda$ that will be  chosen for $\hat r_-$ and for $\hat r_+$
will not depend on $a$.\rm\
 Since $u_{+,a,r}(0)$ and $v_{+,a,r}(0)$ are decreasing in their dependence on $r$,
it will then follow from \eqref{rr+r-} that
\begin{equation}\label{rhatfrac}
\frac{v_{+,a,\hat r_+}(0)}{u_{+,a,\hat r_-}(0)}\le E_0^{(r)}T_a\le\frac{v_{+,a,\hat r_-}(0)}{u_{+,a,\hat r_+}(0)}.
\end{equation}
In the case that $\gamma$ and $\lambda$ have been chosen to construct
$\hat r_-$, denote the function $\phi_3$ above by $\phi_{3,-}$, and in the case that $\gamma$ and $\lambda$ have been chosen to construct
$\hat r_+$, denote the function $\phi_3$ above by $\phi_{3,+}$.
We will then be able to complete the proof of the theorem using \eqref{rhatfrac} along with  Proposition \ref{uv+-}, which gives
 $u_{+,a,\hat r_+}$ and $v_{+,a,\hat r_+}$
 explicitly in terms of $\phi_{3,+}$, and  $u_{+,a,\hat r_-}$ and $v_{+,a,\hat r_-}$
 explicitly in terms of $\phi_{3,-}$.

We begin with finding $\gamma$ and $\lambda$ to construct $\hat r_-$ in the case $l\in(-1,0)$.
This is the most delicate case. The term in the square brackets on the right hand side of \eqref{hatr} will clearly be positive for
all $x$ if $(l+1)\lambda\gamma^{\frac{l+1}2}+l\ge0$; thus, in particular, it will be positive for all $x$ if $\gamma=\gamma(\lambda):=(\frac{-l}{(l+1)\lambda})^{\frac2{l+1}}$. Thus, from \eqref{hatr},
the inequality $\hat r_-\ge0$ will hold for any $\lambda>0$, if we choose $\gamma=\gamma(\lambda)$.
We now show that if $\lambda$ is chosen sufficiently small, and $\gamma=\gamma(\lambda)$, then
$r_-\ge \hat r_-$.  We have
\begin{equation}\label{rminusrhat}
\begin{aligned}
&r_-(x)-\hat r_-(x)=\frac D2\lambda(l+1)(\gamma(\lambda)+x^2)^{\frac{l-3}2}\times\\
&\Big[\frac{2c_1}{D\lambda(l+1)}(\gamma(\lambda)+x^2)^{\frac{3-l}2}(\gamma_1+x^2)^l-
(l+1)\lambda x^2(\gamma(\lambda)+x^2)^{\frac{l+1}2}-\gamma(\lambda)-lx^2\Big].
\end{aligned}
\end{equation}
Thus, it remains to show that for sufficiently small $\lambda$, the expression in the square brackets in \eqref{rminusrhat}
is nonnegative for all $x$. Since for $a,b>0$, one has $(a+b)^m\le a^m+b^m$ if $m\in[0,1]$ and $(a+b)^m\ge a^m+b^m$
if $m\ge 1$,  the expression in the square brackets will be positive if
\begin{equation*}
\frac{2c_1}{D\lambda(l+1)}\big(\gamma_1+x^2\big)^l\big(\gamma(\lambda)^{\frac{3-l}2}+x^{3-l}\big)-(l+1)\lambda x^2\big(\gamma(\lambda)^{\frac{l+1}2}+x^{l+1}\big)-\gamma(\lambda)-lx^2\ge0.
\end{equation*}
Since $\gamma(\lambda)^{\frac{l+1}2}=\frac{-l}{\lambda(l+1)}$, the above inequality can be rewritten as
\begin{equation}\label{analyze}
\frac{2c_1}{D\lambda(l+1)}\big(\gamma_1+x^2\big)^l\big(\gamma(\lambda)^{\frac{3-l}2}+x^{3-l}\big)-(l+1)\lambda x^{3+l}-\gamma(\lambda)\ge0.
\end{equation}
Noting that the first term in \eqref{analyze} behaves asymptotically like
$\frac{2c_1}{D\lambda(l+1)}x^{3+l}$ as $x\to\infty$ and noting that $\gamma(\lambda)\to\infty$ as $\lambda\to0$
and that $\frac{3-l}2>1$, it is easy to see that \eqref{analyze} holds for all $x$.

We now find $\gamma$ and $\lambda$ to construct $\hat r_-$ in the case $l\ge0$.
From \eqref{hatr}, we automatically have $\hat r_-\ge0$.
We fix $\gamma$ arbitrarily and consider small $\lambda$. We have as in \eqref{rminusrhat},
\begin{equation}\label{rminusrhatagain}
\begin{aligned}
&r_-(x)-\hat r_-(x)=\frac D2\lambda(l+1)(\gamma+x^2)^{\frac{l-3}2}\times\\
&\Big[\frac{2c_1}{D\lambda(l+1)}(\gamma+x^2)^{\frac{3-l}2}(\gamma_1+x^2)^l-
(l+1)\lambda x^2(\gamma+x^2)^{\frac{l+1}2}-\gamma-lx^2\Big].
\end{aligned}
\end{equation}
The term in the square brackets on the right hand side of \eqref{rminusrhatagain}, when  evaluated at $x=0$, is equal to
$\frac{2c_1}{D\lambda(l+1)}\gamma^{\frac{3-l}2}\gamma_1^l-\gamma$, and for large $|x|$ behaves asymptotically
like $\big(\frac{2c_1}{D\lambda(l+1)}-(l+1)\lambda\big)x^{l+3}$. From this and the  general form of the term
in the square brackets,
it is clear that the right hand side of \eqref{rminusrhatagain} is positive for all $x$ if $\lambda$ is chosen
sufficiently small.

We now find $\gamma$ and $\lambda$ to construct $\hat r_+$ for any $\l>-1$.
From \eqref{hatr}, one has $\hat r(0)=\frac D2\lambda(l+1)\gamma^{\frac{l-1}2}$ and
$\hat r(x)\sim\frac D2\lambda^2(l+1)^2x^{2l}$ as $x\to\infty$.
It is clear from this and the general form of \eqref{hatr} that if one fixes $\gamma$ arbitrarily and lets $\lambda$ be sufficiently large,
then $\hat r_+(x)\ge r_+(x)$ for all $x$.

We now turn to estimating $u_{+,a,\hat r_\pm}(0)$ and $v_{+,a,\hat r_\pm}(0)$, using
Proposition \ref{uv+-}.
 From
\eqref{u+3} with $\phi_{3,+}$ or $\phi_{3,-}$ in place of $\phi_3$, it is clear that
$u_{+,a,\hat r_\pm}(0)$   satisfy the estimates
\begin{equation}\label{uestimates}
\begin{aligned}
&u_{+,a,\hat r_\pm}(0)\sim \frac{C_\pm}{\phi_{3,+}(a)},\ \text{as}\ a\to\pm \infty,
\ \text{for some}\ C_\pm>0;\\
&\lim_{a\to0^+}u_{+,a,\hat r_\pm}(0)=1.
\end{aligned}
\end{equation}
Now consider \eqref{v+3} with $\phi_{3,+}$ or $\phi_{3,-}$ in place of $\phi_3$.
Since $\int_{-\infty}^{\infty}\phi^{-2}_{3,\pm}(y)dy<\infty$,  the denominator   $\int_{-\infty}^a\phi^{-2}_{3,\pm}(y)dy$ of the fraction
in \eqref{v+3} is bounded as $a\to\infty$.
 Write the numerator of that fraction with $x=0$ as
\begin{equation*}
\begin{aligned}
&\int_{-\infty}^0dy\thinspace\phi^{-2}_{3,\pm}(y)\int_0^adt\phi^{-2}_{3,\pm}(t)\int_y^t\phi_{3,\pm}(z)dz=
\int_{-\infty}^0dy\thinspace\phi^{-2}_{3,\pm}(y)\int_0^adt\phi^{-2}_{3,\pm}(t)\int_0^t\phi_{3,\pm}(z)dz+\\
&\int_{-\infty}^0dy\thinspace\phi^{-2}_{3,\pm}(y)\int_0^adt\phi^{-2}_{3,\pm}(t)\int_y^0\phi_{3,\pm}(z)dz.
\end{aligned}
\end{equation*}
Since the functions $\phi_{3,\pm}$ are increasing on $[0,\infty)$ and decreasing on $(-\infty,0]$, we  have
$$
\begin{aligned}
&\int_{-\infty}^0dy\thinspace\phi^{-2}_{3,\pm}(y)\int_0^adt\phi^{-2}_{3,\pm}(t)\int_0^t\phi_{3,\pm}(z)dz
\le\big(\int_{-\infty}^0\phi^{-2}_{3,\pm}(y)dy\big)\big(\int_0^\infty  t\phi^{-1}_{3,\pm}(t)dt\big)<\infty;\\
&\int_{-\infty}^0dy\thinspace\phi^{-2}_{3,\pm}(y)\int_0^adt\phi^{-2}_{3,\pm}(t)\int_y^0\phi_{3,\pm}(z)dz\le
\big(\int_{-\infty}^0y\phi^{-1}_{3,\pm}(y)dy\big)\big(\int_0^\infty  \phi^{-2}_{3,\pm}(t)dt\big)<\infty.
\end{aligned}
$$
We conclude from this that $v_{+,a,\hat r_\pm}(0)$ are bounded as $a\to\infty$. It is also clear from \eqref{v+3} that
$v_{+,a,\hat r_\pm}(0)\sim ca$ as $a\to0^+$, for some $c>0$.
Using these facts  along  with
\eqref{uestimates}, \eqref{rhatfrac} and \eqref{phis}, we conclude that \eqref{solugrowth} holds.

\hfill $\square$

\section{Proof of Theorem  \ref{comparison}}\label{proofthmcomp}
We will prove the theorem for $a>0$; the same type of proof works for $a<0$.
As in the proof of Theorem \ref{comparison1},
since we will need to compare with different choices of $r$,
we denote the
solutions $u_{+,a}$ and $v_{+,a}$ of \eqref{u+a}  and \eqref{v+a}
 by $u_{+,a,r}$ and $v_{+,a,r}$.

\noindent \it Parts  (i) and (ii).\rm\
 By Proposition \ref{voveru}, $E_0^{(r)}T_a=\infty$ if and only  if
$v_{+,a}(0)=\infty$. Define
$$
\phi_3(x)=\gamma+x^2,\ \gamma>0.
$$
 Let
$$
\hat r(x):=\frac D2\frac{\phi_3''(x)}{\phi_3(x)}=\frac D{\gamma+x^2}.
$$
Then  $v_{+,a,\hat r}$ is given as in \eqref{v+3} with $\phi_3$
as above. The numerator of the fraction in \eqref{v+3} with $x=0$ satisfies
$$
\begin{aligned}
&\int_{-\infty}^0dy\thinspace\phi^{-2}_3(y)\int_0^adt\phi^{-2}_3(t)\int_y^t\phi_3(z)dz\ge\\
&\big(\int_{-\infty}^0dy\thinspace\phi^{-2}_3(y)\int_y^0\phi_3(z)dz\big)\big(\int_0^adt\phi^{-2}_3(t)\big)=\infty;\\
\end{aligned}
$$
thus, $v_{+,a,\hat r}(0)=\infty$. Since $v_{+,a,r}(0)$ is decreasing in its dependence on $r$, we conclude
that $E_0^{(r)}T_a=\infty$, if $r(x)\le \frac D{\gamma+x^2}$, for some $\gamma>0$ and for all $x$.

Now $v_{+,a,r}$ can be written as $v_{+,a,r}(x)=\int_{-\infty}^aG_{a,r}(x,y)dy$,
where $G_{a,r}(x,y)$ is the Green's function for the operator $\frac D2\frac{d^2}{dx^2}-r(x)$ on $(-\infty,a)$
\cite{P}. If $r_1,r_2\gneqq0$ and $r_1-r_2$ is compactly supported, then
there exists a constant $c\in(0,1)$ such that
$c\le \frac{G_{a,r_1}(x,y)}{G_{a,r_2}(x,y)}\le \frac1c$, for all $x,y\in(-\infty, a)$ \cite{Pinch}.
It follows from this that the finiteness or infiniteness of $v_{+,a,r}$ is not affected by
compactly supported changes in $r$.
Thus, we conclude that $E_0^{(r)}T_a=\infty$, if $r(x)\le \frac D{\gamma+x^2}$,
for some $\gamma>0$ and
for all sufficiently large
$|x|$. This proves part (i).

We turn to part (ii). Let $\lambda>1$ and $\gamma>0$.
Define $\phi_3(x)=\gamma_1+|x|^m$, where
$\gamma_1>0$ and
$m>2$ is chosen
so that $\lambda_0:=\frac{m(m-1)}2$ satisfies $\lambda_0<\lambda$. Let
$$
\hat r(x):=\frac D2\frac{\phi_3''(x)}{\phi_3(x)}=\frac{D\lambda_0|x|^{m-2}}{\gamma_1+|x|^m}.
$$
Then $v_{+,a,\hat r}$ is given  as in \eqref{v+3} with $\phi_3$ as above. This
shows that $v_{+,a,\hat r}(0)<\infty$.
 If $\gamma_1$ is chosen sufficiently large, then $\frac{D\lambda}{\gamma+x^2}\ge\hat r(x)$.
Thus $v_{+,a,\frac{D\lambda}{\gamma+x^2}}(0)<\infty$, since $v_{+,a,r}$ is decreasing in its dependence
on $r$. Since the finiteness or infiniteness of $v_{+,a,r}$ is not affected by
compactly supported changes in $r$, we conclude that $v_{+,a,r}(0)<\infty$ if
$r(x)\ge \frac {D\lambda}{\gamma+x^2}$ for sufficiently large $|x|$, and thus also $E_0^{(r)}T_a<\infty$ for such $r$.
\medskip

\noindent \it Part (iii).\rm\    Fix $\lambda>1$ and choose
$m=\frac{1+\sqrt{1+8\lambda}}2>2$ so that   $\frac{m(m-1)}2=\lambda$. Define
$\phi_3(x)=1+|x|^m$ and let $r(x):=\frac D2\frac{\phi''_3(x)}{\phi_3(x)}=\frac{D\lambda|x|^{m-2}}{1+|x|^m}$.
Then $r(x)\sim \frac{D\lambda}{|x|^2}$ as $|x|\to\infty$.
Using  \eqref{u+3} and \eqref{v+3} with $\phi_3$ as above, it follows that $u_{+,a,r}(0)\sim \frac {C_0}{\phi_3(a)}$ as $a\to\infty$, for some $C_0>0$, and that $v_{+,a,r}$ is bounded as $a\to\infty$. Thus, it follows from  Proposition
\ref{voveru} that
$E_0^{(r)}T_a\sim C|a|^{\frac{1+\sqrt{1+8\lambda}}2}$, for some $C>0$.
\medskip

\it\noindent Part (iv).\rm\ Fix $\lambda_1,\lambda_2,\gamma_1,\gamma_2$ as in the statement of the theorem.
Fix  $\epsilon>0$, and define $m=\frac{1+\sqrt{1+8\lambda_1}}2-\epsilon$
and $M=\frac{1+\sqrt{1+8\lambda_2}}2+\epsilon$. Assume that $\epsilon$ is sufficiently small so that
$m>2$.
We will show that  $c_1,c_2>0$ can be chosen  so that $\phi_{3,-}(x):=c_1+c_2x^2+|x|^{m}$
satisfies
\begin{equation}\label{r-c1c2}
r_-(x):=\frac D2\frac{\phi_{3,-}''(x)}{\phi_{3,-}(x)}\le \frac{D\lambda_1}{\gamma_1+x^2},
\end{equation}
and that $c_1,c_2>0$ can be chosen so that
$\phi_{3,+}(x):=c_1+c_2x^2+|x|^{M}$
satisfies
\begin{equation}\label{r+c1c2}
r_+(x):=\frac D2\frac{\phi_{3,+}''(x)}{\phi_{3,+}(x)}\ge \frac{D\lambda_2}{\gamma_2+x^2}.
\end{equation}
It  will then follow from Proposition \ref{voveru} and the fact that $u_{+,a,r}$ and $v_{+,a,r}$ are decreasing
in their dependence on $r$ that
\begin{equation}\label{r+-frac}
\frac{v_{+,a,r_+}(0)}{u_{+,a,r_-}(0)}\le E_0^{(r)}T_a\le \frac{v_{+,a,r_-}(0)}{u_{+,a,r_+}(0)}
\end{equation}
The functions $u_{+,a,r_\pm}$ are given by \eqref{u+3} with $\phi_3$ replaced by $\phi_{3,\pm}$ from above,
and the functions $v_{+,a,r_\pm}$ are given by \eqref{v+3} with $\phi_3$ replaced by $\phi_{3,\pm}$.
One finds that $v_{+,a,r_\pm}(0)$ are bounded as $a\to\infty$ and that
$u_{+,a,r_\pm}(0)\sim \frac{C_\pm}{\phi_{3,\pm}(a)}$
as $a\to\infty$, for constants $C_\pm>0$. Using this with \eqref{r+-frac}  proves \eqref{witheps}.

It remains to find a pair $c_1,c_2$ for $\phi_{3,-}$ and a pair $c_1,c_2$ for $\phi_{3,+}$.
Define $\delta_1,\delta_2>0$ by
$\frac{m(m-1)}2=\lambda_1-\delta_1$ and
$\frac{M(M-1)}2=\lambda_2+\delta_2$.
We begin with $\phi_{3,-}$. We have
$$
\frac D2\frac{\phi_{3,-}''}{\phi_{3,-}}=\frac{Dc_2+D(\lambda_1-\delta_1)|x|^{m-2}}{c_1+c_2x^2+|x|^{m}}.
$$
Thus, the inequality \eqref{r-c1c2} we wish to satisfy  can be written as
\begin{equation}\label{r-c1c2need}
\begin{aligned}
&D\lambda_1c_1+D\lambda_1c_2x^2+D\lambda_1|x|^{m}\ge\\ &Dc_2\gamma_1+Dc_2x^2+D(\lambda_1-\delta_1)\gamma_1|x|^{m-2}+D(\lambda_1-\delta_1)|x|^{m}.
\end{aligned}
\end{equation}
It is clear that we can choose $c_2$ sufficiently large so that
$D\lambda_1c_2x^2+D\lambda_1|x|^{m}\ge Dc_2x^2+D(\lambda_1-\delta_1)\gamma_1|x|^{m-2}+D(\lambda_1-\delta_1)|x|^{m}$.
Once such a $c_2$ is chosen, it is clear that $c_1$ can be chosen sufficiently large so that \eqref{r-c1c2need} holds.

We now find a pair $c_1,c_2$ for $\phi_{3,+}$. We have
$$
\frac D2\frac{\phi_{3,+}''}{\phi_{3,+}}=\frac{Dc_2+D(\lambda_2+\delta_2)|x|^{M-2}}{c_1+c_2x^2+|x|^{M}}.
$$
Thus, the inequality \eqref{r+c1c2} we wish to satisfy  can be written as
\begin{equation}\label{r+c1c2need}
\begin{aligned}
&D\lambda_2c_1+D\lambda_2c_2x^2+D\lambda_2|x|^M\le\\ &Dc_2\gamma_2+Dc_2x^2+D(\lambda_2+\delta_2)\gamma_2|x|^{M-2}+D(\lambda_2+\delta_2)|x|^M.
\end{aligned}
\end{equation}
It is clear that we can choose $c_2$ sufficiently small so that
$$
\begin{aligned}
&\inf_{x\in\mathbb{R}}\Big(Dc_2\gamma_2+Dc_2x^2+D(\lambda_2+\delta_2)\gamma_2|x|^{M-2}+D(\lambda_2+\delta_2)|x|^M-\\
&D\lambda_2c_2x^2-D\lambda_2|x|^M\Big)>0.
\end{aligned}
$$
Once $c_2$ has been chosen, it is clear that $c_1$ can be chosen sufficiently small so that \eqref{r+c1c2need}
holds.
\hfill $\square$

\section{When $\mu$ is supported on $[-L_1,L_2]$}\label{finsupp}
Let  $L_1,L_2>0$ and assume that the support of the target distribution is $[-L_1,L_2]$. In this case, there is no reason to search outside of the above interval, and thus as soon as the searcher reaches
$-L_1$ or $L_2$, its position should be reset to 0. This is equivalent to setting $r\equiv\infty$ off of $[-L_1,L_2]$.
On  $[-L_1,L_2]$ we only need assume now
that $r\ge0$, not that $r\gneqq0$.
We'll use the notation $P_x^{(r;[-L_1,L_2])}, E_x^{(r;[-L_1,L_2])}$ for probabilities and expectations for this process starting from $x\in[-L_1,L_2]$.
As in section \ref{proofpropvoveru}, for fixed $a\in(0,L_2]$ let
 $\mathcal{T}_t$ be the semigroup defined by
$$
\begin{aligned}
&\mathcal{T}_tf(x)=E_x^{(r;[-L_1,L_2])}(f(X(t));T_a>t):=E_x^{(r;[-L_1,L_2])}(f(X(t))1_{\{T_a>t\}}),\\
& x\in(-L_1,a],
\end{aligned}
$$
for bounded continuous $f$, with a parallel  definition in the case $a\in[-L_1,0)$.
Its generator is $\mathcal{L}$ as in \eqref{generator} with the zero Dirichlet boundary condition at $x=a$ and with the additional boundary condition requiring that the value of the function
at $-L_1$ be equal to the value of the function at 0. (For boundary conditions for semigroups corresponding to processes that jump from the boundary to an interior point, see for example \cite{BAP}.)

Let $w(x,t)=\mathcal{T}_t1(x)$. Then $w(x,t)=P_x^{(r;[-L_1,L_2])}(T_a>t)$ and it solves
\begin{equation}\label{w}
\begin{aligned}
&w_t=\mathcal{L}w=\frac D2w''+r(x)\big(w(0,t)-w(x,t)\big), \ x\in(-L_1,a);\\
&w(x,0)=1,\ x\in(-L_1,a);\\
&w(a,t)=0,\   w(-L_1,t)=w(0,t),\  t>0.
\end{aligned}
\end{equation}
In the present case, the statement of the  result analogous to   Proposition \ref{voveru} from section \ref{proofpropvoveru} looks exactly the same, the only difference being that
the functions $u_{\pm,a}$ and $v_{\pm,a}$
that satisfied the
equations \eqref{u+a}-\eqref{v-a} will now be called $u_{\pm,a;-L_1,L_2}$ and $v_{\pm,a;-L_1,L_2}$ and
they
 satisfy the  equations below instead, \eqref{u+aL}-\eqref{v-aL}. The only difference in the proof is that since the interval we work with is bounded---$(-L_1,a)$, if $a>0$ and  $(a,L_2)$ if $a<0$, we
 don't need to truncate to a bounded domain as was done in the proof of Proposition \ref{voveru}.

\noindent For $a\in(0,L_2]$, let $u_{+,a;-L_1,L_2}$ denote the solution to the equation
\begin{equation}\label{u+aL}
\begin{aligned}
&\frac D2 u''-r(x)u=0, \ x\in(-L_1,a);\\
&u(a)=1;\\
&u(-L_1)=u(0).
\end{aligned}
\end{equation}
For $a\in[-L_1,0)$,  let $u_{-,a;-L_1,L_2}$ denote the solution to
\begin{equation}\label{u-aL}
\begin{aligned}
&\frac D2 u''-r(x)u=0, \ x\in(a,L_2);\\
&u(a)=1;\\
&u(L_2)=u(0).
\end{aligned}
\end{equation}
For $a\in(0,L_2]$, let $v_{+,a;-L_1,L_2}$ denote the solution to
\begin{equation}\label{v+aL}
\begin{aligned}
&\frac D2 v''-r(x)v=-1,\  x\in(-L_1,a)  ;\\
&v(a)=0;\\
&v(-L_1)=v(0).
\end{aligned}
\end{equation}
For $a\in[-L_1,0)$  let $v_{-,a;-L_1,L_2}$ denote the solution to
\begin{equation}\label{v-aL}
\begin{aligned}
&\frac D2 v''-r(x)v=-1,\ x\in(a,L_2) ;\\
&v(a)=0;\\
&v(L_2)=v(0).
\end{aligned}
\end{equation}
We record the result that corresponds to Proposition \ref{voveru}.
\begin{proposition}\label{voveruinfinity}
Let $r\ge0$ be a continuous function on $[-L_1,L_2]$. Then
\begin{equation}\label{voveruformula}
E_0^{(r;[-L_1,L_2])}T_a=\begin{cases}\frac{v_{+,a;-L_1,L_2}(0)}{u_{+,a;-L_1,L_2}(0)}, \ 0<a\le L_2;\\
\frac{v_{-,a;-L_1,L_2}(0)}{u_{-,a;-L_1,L_2}(0)}, \ -L_1\le a<0.\end{cases}
\end{equation}
\end{proposition}
Similar to Proposition \ref{criticality}, for any $r\ge0$, one can find a positive solution $\phi$ to $\frac D2\Delta\phi-r(x)\phi=0$ in $[-L_1,L_2]$.
Using such a function $\phi$ in the manner that we used $\phi_i$, $i=1,2,3$, one can proof a result parallel to Proposition \ref{uv+-} that
gives the solutions  $u_{\pm,a;-L_1,L_2}$ and $v_{\pm,a;-L_1,L_2}$ explicitly in terms of $\phi$.
As in Theorem \ref{rep}, one then obtains an explicit formula for $E_0^{(r;[-L_1,L_2])}T_a$ in terms
of this function $\phi$.
The formulas are a bit more complicated than those appearing
in Proposition \ref{uv+-}, so we refrain from writing them down.

However, we will consider in detail the case that $r(x)=r\ge0$ is constant on $[-L_1,L_2]$.
Assume first that $r>0$.
In this case it is easy to see from the equations that $v_{+,a;-L_1,L_2}=\frac1r(1-u_{+,a;-L_1,L_2})$.
Thus, from \eqref{voveruformula},
$E_0^{(r;[-L_1,L_2])}T_a=\frac1r(\frac1{u_{+,a;-L_1,L_2}(0)}-1)$, for $a\in(0,L_2]$.
Solving for $u_{+,a;-L_1,L_2}$ from \eqref{u+aL} and substituting $x=0$, and making  similar calculations for $u_{-,a;-L_1,L_2}$,
we obtain
\begin{equation}\label{constantrinterval}
E_0^{(r;[-L_1,L_2])}T_a=\begin{cases}
\frac1r\Big[\frac{\sinh\sqrt{\frac{2r}D}(a+L_1)-\sinh\sqrt{\frac{2r}D}a-\sinh\sqrt{\frac{2r}D}L_1}{\sinh\sqrt{\frac{2r}D}L_1}\Big],\ 0<a\le L_2\\
\frac1r\Big[\frac{\sinh\sqrt{\frac{2r}D}(|a|+L_2)-\sinh\sqrt{\frac{2r}D}|a|-\sinh\sqrt{\frac{2r}D}L_2}{\sinh\sqrt{\frac{2r}D}L_2}\Big],\ -L_1\le a<0.
\end{cases}
\end{equation}
Making similar calculations when $r=0$, or simply taking the limit of the above expression when $r\to0$, one obtains
\begin{equation}\label{r=0}
E_0^{(r;[-L_1,L_2])}T_a=\begin{cases}\frac{a(a+L_1)}D,\ 0<a\le L_2;\\
\frac{|a|(|a|+L_2)}D,\ -L_1\le a<0.\end{cases}
\end{equation}

From now on, consider the symmetric case, $L_1=L_2=A$, for some $A>0$.
Consider the uniform target distribution   on  $[-A,A]$.  The expected  distance to the target
is then AvgDist$:=\frac A2$.
The expected time to locate the target is
$$
\frac1{2A}\int_{-A}^AE_0^{(r;[-A,A])}T_a\thinspace da=\frac1A\int_0^AE_0^{(r;[-A,A])}T_a\thinspace da.
$$
Substituting from \eqref{constantrinterval} and performing the integration, we obtain
$$
\begin{aligned}
&\frac1{2A}\int_{-A}^AE_0^{(r;[-A,A])}T_a\thinspace da=\\
&\frac1{rA\sinh(\sqrt{\frac{2r}D}A)}\sqrt{\frac D{2r}}\Big(\cosh(\sqrt{\frac{2r}D}(2A)-2\cosh(\sqrt{\frac{2r}D}A)+1\Big)-\frac1r.
\end{aligned}
$$
Letting $x=\sqrt{\frac{2r}D}A$, we have
$$
\frac1{2A}\int_{-A}^AE_0^{(r;[-A,A])}T_a\thinspace da=\frac{2A^2}D\Big(\frac{\cosh(2x)-2\cosh(x)+1}{x^3\sinh(x)}-\frac1{x^2}\Big).
$$
The minimum of the function in parentheses above is obtained at $x=0$, the value of the function there being
$\frac5{12}$. Thus we conclude that the expected  time to locate the target is minimized in the class
of constant resetting rates $r$ on $[-A,A]$  by setting $r=0$, and
$$
\inf_{r\ge0, r\ \text{constant}}\frac1{2A}\int_{-A}^AE_0^{(r;[-A,A])}T_a\thinspace da=\frac56\frac{A^2}{D}=
\frac{10}3\frac{(\text{AvgDist})^2}D.
$$

Intuitively, it seems then that the minimum will also be obtained at $r=0$ if the symmetric
target distribution has a non-decreasing density on $[0,A]$. We now consider
the linearly decreasing, symmetric density which decreases to zero.
The expected distance to the target is then
$\frac2{A^2}\int_0^Aa(A-a)da=\frac13A$.
We have
$$
\int_{-A}^A(E_0^{(r;[-A,A])}T_a)\thinspace \frac{A-|a|}{A^2}da=\frac2{A^2}\int_0^A(E_0^{(r;[-A,A])}T_a)\thinspace (A-a)da.
$$
Substituting from \eqref{constantrinterval} and performing the integration, and again making the substitution
$x=\sqrt{\frac{2r}D}A$, we obtain
$$
\int_{-A}^A(E_0^{(r;[-A,A])}T_a)\thinspace \frac{A-|a|}{A^2}da=\frac{4A^2}D\Big(\frac{\frac{\sinh(2x)}x-\frac{2\sinh(x)}x-\cosh(x)+1}{x^3\sinh(x)}-\frac1{x^2}\Big).
$$
The minimum of the function in the parentheses above is obtained at $x\approx1.3538$ and the minimum value
is approximately $0.1238$.
Thus we conclude that the expected  time to locate the target is minimized in the class
of constant resetting rates $r$ on $[-A,A]$  by setting $r\approx0.916\frac D{A^2}$, and
$$
\inf_{r\ge0, r\ \text{constant}}
\int_{-A}^A(E_0^{(r;[-A,A])}T_a)\thinspace \frac{A-|a|}{A^2}da\approx0.495\frac{A^2}D=4.455\frac{(\text{AvgDist})^2}D.
$$

It might be interesting to pursue the above direction of calculations further. In particular, what can
be said about the ratio
$$
\begin{aligned}
&\frac{D\inf_{r\ge0, r\ \text{constant}}
\int_{-A}^A(E_0^{(r;[-A,A])}T_a)\thinspace \mu(da)}{(\text{AvgDist}(\mu))^2}=\\
&\frac{D\inf_{r>0}
\int_{-A}^A\frac1r\Big[\frac{\sinh\sqrt{\frac{2r}D}(a+A)-\sinh\sqrt{\frac{2r}D}a-\sinh\sqrt{\frac{2r}D}A}{\sinh\sqrt{\frac{2r}D}A}\Big] \mu(da)}{(2\int_0^Aa\mu(da))^2},
\end{aligned}
$$
as one varies over all symmetric distributions $\mu$ with support  $[-A,A]$?

\bf\noindent Acknowledgment.\rm\ The author thanks his colleague, Nir Gavish, for the numerical analysis on page 9.

\end{document}